\documentclass[a4paper,12pt]{article}

\usepackage{amsmath,amsfonts,latexsym,theorem}
\nonstopmode

\numberwithin{equation}{section}

\newtheorem{Theorem}{Theorem}[section]
\newtheorem{Definition}{Definition}[section]
\newtheorem{Proposition}{Proposition}[section]
\newtheorem{Lemma}{Lemma}[section]

\newenvironment{Proofc}[1]{\smallskip\par\noindent\textsc{#1}\quad}%
  {\hfill$\Box$\bigskip\par}
\newenvironment{Proof}{\begin{Proofc}{Proof}}{\end{Proofc}}

\newtheorem{Remark}{Remark}[section]

\def\a{\alpha}
\def\b{\beta}
\def\d{\delta}

\def\g{\gamma}
\def\G{\Gamma}
\def\l{\lambda}

\def\m{\mu}
\def\n{\nu}

\def\th{\theta}

\def\t{\tau}

\def\pd{\partial}

\newcommand{\cV}{{\cal V}}
\newcommand{\cE}{{\cal E}}

\newcommand{\ds}{\displaystyle}
\newcommand{\Gin}{{\cV}}

\newcommand{\cK}{{\cal K}}

\newcommand{\R}{{\mathbb R}}
\newcommand{\N}{{\mathbb N}}

\def\pd{\partial}

\begin{document}
\title{Stationary  Mean Field Games systems defined  on   networks}
 \author{Fabio Camilli\footnotemark[1] \and Claudio Marchi\footnotemark[2]}

\footnotetext[1]{Dip. di Scienze di Base e Applicate per l'Ingegneria,  ``Sapienza" Universit{\`a}  di Roma, via Scarpa 16,
 00161 Roma, Italy, ({\tt e-mail:camilli@dmmm.uniroma1.it})}
\footnotetext[2]{Dip. di Ingegneria dell'Informazione, Universit\`a di Padova, via Gradenigo 6/B, 35131 Padova, Italy ({\tt claudio.marchi@unipd.it}).}

\date{version: \today}
\maketitle
\begin{abstract}
We consider a stationary Mean Field Games system defined on a network.
In this framework, the transition conditions at the vertices play  a crucial role:
the ones here considered  are based on the optimal control interpretation of the problem.
We prove separately the well-posedness for each of the two equations composing the system.
Finally, we prove  existence and uniqueness of the solution of the Mean Field Games system.
\end{abstract}
  \begin{description}
\item [\textbf{MSC 2000}: ]35R02, 35B30, 49N70.
  \item [\textbf{ Keywords}:] mean field games, networks,  transition conditions.
  \end{description}
%
%
%
\section{Introduction}\label{intro}

The theory of  Mean Field Games (briefly, MFG) has been introduced  in  \cite{hmc, ll} to  describe the asymptotic behavior of  stochastic differential game problems (Nash equilibria) as the number of players tends to $+\infty$.
>From a mathematical  point of view, MFG theory leads to the study of a coupled system of two differential equations: one  equation is of Hamilton-Jacobi-Bellman (briefly, HJB) type and it describes the optimal behavior of the single agent, while  the other one is  a Fokker-Planck (briefly, FP) equation governing the   distribution of the overall population. The system can be completed with different boundary conditions (periodic, Dirichlet, Neumann) and initial conditions (initial-terminal condition, planning problem). Existence and  uniqueness of strong and
weak solutions to the MFG system have  been obtained under rather general assumptions on the data of the problems (\cite{cgr,cgpt,gs,p}).\par
Aim of this note is to study MFG  systems defined on a  network.   While the differential equations are defined in the usual way along  the edges, a crucial issue is to find the correct conditions  at the vertices (\emph{transition conditions}). We will choose a set of transition conditions according to the optimal control interpretation of the system.

Starting with the seminal paper by Lumer \cite{lu}, a general   theory for linear and semilinear differential  equations on networks has been developed mainly employing  the variational structure of the problem. In this framework the natural transition conditions are, besides the continuity of the solution, the so-called Kirchhoff conditions on the first order derivatives (see \cite{mug,nic,pb,vb}).\par
For a single nonlinear equation, existence and uniqueness results  are available  only for some specific classes of operators such as conservation law \cite{cg} and some Hamilton-Jacobi equations \cite{cms}. Therefore the  first  step in our analysis  is to establish existence and uniqueness results for each of the two equations composing the MFG system on the network. Because of its control theoretic interpretation (see \cite{fs,fw} and Section \ref{sec2}) the natural transition condition for the HJB equation is the Kirchhoff condition, while for the FP equation it is natural to require the conservation of the flux at the vertices. In Section \ref{sec2} we discuss the relationship between the transition conditions for the two equations.\par
After having solved the   two   equations separately, we tackle our second and main issue: the study of  the MFG system on a network. We shall obtain existence of a solution by a fixed point argument; moreover we shall get uniqueness of the solution adapting a classical argument to this framework and taking advantage of the  relation between the transition conditions of the two equations.\par
As far as we know, this is the first paper to consider  general MFG systems (HJB equation and FP equation as well) on networks. Indeed, in \cite{ccm} only a particular class of MFG systems on networks was addressed; in that setting, by a suitable change of variables, the HJB and the FP equation are transformed in two heat equations both with Kirchhoff transition conditions coupled via the initial data.
Moreover, it is worth to observe that the papers \cite{g.ppt}, \cite{gms1} and \cite{gms2}  consider MFG systems on graphs (namely,   the state variable belongs to a discrete set).\par
We remark that the results here contained can be generalized in several directions (nonlocal coupling, evolutive problems, boundary conditions, weak solutions, ramified spaces etc); moreover the assumptions are far from  being optimal. Since this is a  first approach to the study of MFG systems on networks, we have tried to keep the presentation as simple as possible in order to avoid technical complications and to concentrate  on the network aspects of the problems.\par
The paper is organized as follows. In the rest of this section, we shall introduce the definition of networks.
In  Section \ref{sec2} we  give a formal derivation of MFG systems on networks and, especially, of the transition conditions. Section \ref{sec3} is devoted to our main results for the HJB equation, the FP equation and, mainly, the MFG system. 
Finally, in the Appendix \ref{appendix} we collect some technical  results.

\vskip 8mm
{\bf Networks.} The network  $\G=(\cV,\cE)$     is a finite  collection of points $\cV:=\{v_i\}_{i\in I}$ in $\R^n$ connected by continuous, non self-intersecting edges $\cE:=\{e_j\}_{j\in J}$. Each edge $e_j\in \cE$ is  parametrized by a smooth function $\pi_j:[0,l_j]\to\R^n,\, l_j>0$. Given $v_i\in \cV$,  we denote by
$Inc_i:=\{j\in J:v_i\in e_j \}$  the set of edges branching out from $v_i$ and by $d_{v_i}:=|Inc_i|$ the degree of $v_i$. A vertex $v_i$ is said a  boundary vertex if $d_{v_i}=1$, otherwise  it is said a transition vertex. For simplicity, in this paper we will assume that the set of boundary vertices is empty.

For a function $u:  \G\to\R$  we denote by $u_j:[0,l_j]\to \R$ the restriction of $u$ to $e_j$, i.e. $u(x)=u_j(y)$ for $x\in e_j$, $y=\pi_j^{-1}(x)$, and  by  $\pd_j u(v_i)$  the  oriented derivative of $u$ at $v_i$ along the arc $e_j$ defined   by
\[
\pd_j u (v_i)=\left\{
               \begin{array}{ll}
                 \lim_{h\to 0^+}(u_j(h)-u_j(0))/h, & \hbox{if $v_i=\pi_j(0)$;} \\
                 \lim_{h\to 0^+}(u_j(l_j-h)-u_j(l_j))/h, & \hbox{if $v_i=\pi_j(l_j)$.}
               \end{array}
             \right.
\] The integral of a function $u$ on $\G$ is defined by
\[\int_\G u(x)dx:=\sum_{j\in J} \int_0^{l_j}u_j(r)dr.\]
The space $L^p(\G)$, $p\ge 1$, is the set    of functions $u:\G\to\R$ such that $u_j\in L^p(0,l_j)$ for all $j\in J$ and  $\|u\|_{p}:=\sum_{j\in J}\|u_j\|_{L^p(0,l_j)}<\infty$. For $p\ge 1$ and for an integer $m>0$, we define the Sobolev space $W^{m,p}(\G)$ as the space of continuous functions on $\G$ such that $u_j\in W^{m,p}(0,l_j)$ for all $j\in J$ and $\|u\|_{m,p}:=\sum_{j\in J}\|u_j\|_{W^{m,p}(0,l_j)}<\infty$.
As usual we set $H^k(\G):=W^{k,2}(\G)$, $k\in\N$.
The space  $C^{k}(\G)$, $k\in \N$, consists of all the continuous functions $u:\G\to\R$  such that  $u_j\in C^{k}([0,l_j])$ for $j\in J$ and    $\|u\|_{C^k}=\max_{\b\le k}\|\pd^\beta u\|_{L^\infty  }<\infty$. Observe that no continuity condition at the     vertices  is   prescribed for the derivatives neither for a function $u\in W^{m,p}(\G)$ nor for a  function $u\in C^k(\G)$.\\
Finally the space  $C^{k,\a}(\G)$, for  $k\in \N$ and $\a\in (0,1)$, is the space of functions $u\in C^k(\G)$ such that $\pd^ku_j\in C^{0,\a}([0,l_j])$ for any $j\in J$ with the norm
\[
\|u\|_{C^{k,\a}}:=\|u\|_{C^k}+\sup_{j\in J}\sup_{x,y\in[0,l_j]}[|\pd^k u_j(x)-\pd^k u_j(y)|/|x-y|^\a].\]

\section{A formal derivation of the MFG system}\label{sec2}

The MFG system  can be deduced from two different points of view (see \cite{ll}): either as the characterization of a  Pareto  equilibrium for dynamic  games with a large number of (indistinguishable) players; or as the optimality conditions for an optimal control problem whose dynamic is governed by a PDE. We explain the two different points of view for MFG systems on networks showing that they lead  to the same   transition conditions. \par

\vskip 8pt
\emph{Pareto equilibrium:}
Consider a population of indistinguishable agents, distributed at time $t=0$ according to the probability $m_0$; any agent moves on a network $\G$ and its dynamics  inside the edge $e_j$ is governed by  the stochastic differential equation
 \[
 dX_s = -\gamma_s\,ds+\sqrt{2\nu_j}\,dW_s,
 \]
 where $\gamma$ is the control, $\nu_j>0$ and  $W_t$ is a 1-dimensional Brownian motion.
When the agent reaches a vertex $v_i\in \Gin$, it almost surely spends zero time at $v_i$ and  enters in one of the incident edges, say $e_j$ with $j\in Inc_i$, with probability $\beta_{ij}$ where
\[\beta_{ij}>0,\, \sum_{j\in Inc_i}\beta_{ij}=1.\]
(see \cite{fs, fw} for a rigorous definition of  stochastic  processes on networks). The cost criterion is given by
\[
 \mathbb{E}_{x}\left[\int_0^T\{L(X_t,\g_t)+V[m(X_t)]\}dt +V_0[m(X_T)]\right]
\]
where $m$ represents the distribution of the overall population of players.
A formal application of the dynamic programming principle gives that the value function $u$ of the previous control problem
satisfies
\begin{equation}\label{MFG1}
\left\{
\begin{array}{ll}
-u_t-\nu_j \pd^2 u  +H_j(x, \pd u)=V[m],\qquad &(x,t)\in e_j\times (0,T),\,j\in J\\[4pt]
\sum_{j\in Inc_i}\a_{ij}\nu_j\pd_j u(v_i,t)=0\quad &(v_i,t)\in\Gin\times (0,T),\\[4pt]
u_j(v_i,t)=u_k(v_i,t),&j,k\in Inc_i,\, (v_i,t)\in\Gin\times(0,T),\\[4pt]
u(x,T)=V_0[m(T)], \quad &x\in\G
\end{array}
\right.
\end{equation}
where $\a_{ij}:=\beta_{ij}\nu_j^{-1}$ and the Hamiltonian is given on the edge $e_j$ by
\[
H_j(x,p)= \sup_{\gamma} \,\big[ -\gamma \cdot p  -L_j(x,\gamma)\big].
\]
Note that the differential equation  inside $e_j$ is defined in terms of  the coordinate parametrizing  the edge.
The second equation in \eqref{MFG1} is known as the Kirchhoff transition condition and it is consequence of the assumption on the  behavior of $X_t$ at the vertices (see \cite{fw}). The third line amounts to the continuity at transition vertices.\par
In order to derive the equation satisfied by the distribution $m$ of the agents, we follow a duality argument. 
Consider the linearized Hamilton-Jacobi  equation
\begin{equation}\label{hjblin}
\left\{
\begin{array}{ll}
    - w_t-\nu \pd^2 w  +\pd_p H(x, \pd u)\pd w=0,\qquad &(x,t)\in e_j\times (0,T),\,j\in J\\
    \sum_{j\in Inc_i}\nu_j\a_{ij}\pd_jw(v_i,t)=0  & (v_i,t)\in\Gin\times(0,T)\\
    w_j(v_i,t)=w_k(v_i,t),&j,k\in Inc_i,\, (v_i,t)\in\Gin\times(0,T),\\
    w(x,T)=0&x\in\G
\end{array}
\right.
\end{equation}
Writing the weak formulation of  \eqref{hjblin} for a test function $m$,  integrating by parts along each edge and regrouping  the boundary terms corresponding to the same vertex $v_i$, we get
\begin{align*}
    0=&\sum_{j\in J} \int_0^T\int_{e_j} \big(-w_t-\nu_j \pd^2 w+\pd_p H_j(x,\pd u)\pd w\big)m\,dxdt\\
      =&\sum_{j\in J}\Big( \int_{e_j}[wm]_0^Tdx+\int_0^T\int_{e_j} \big[ m_t -\nu_j \pd^2 m -\pd (m\pd_p H_j(x,\pd u))\big]w\,dxdt\Big)\\
      -&\int_0^T\sum_{v_i\in\Gin}\left[\sum_{j\in Inc_i}\nu_j m_j(v_i,t)\pd_j w (v_i,t)-\big( \nu_j \pd_j m(v_i,t)+ \pd_pH(v_i,\pd u)m_j(v_i,t)\big)w (v_i,t)\right]dt.
\end{align*}
By the previous identity   we obtain that $m$ satisfies inside the edges the  equation
\[
  m_t -\nu \pd^2 m -\pd(m\,\pd_p H(x, \pd u))=0.
 \]
Moreover, recalling  the transition condition for $w$,  the  first one   of the three terms computed at the transition vertices    vanishes if
\begin{equation}\label{condm1}
  \frac{m_j(v_i,t)}{\a_{ij}}= \frac{m_k(v_i,t)}{\a_{ik}}, \qquad j,k\in Inc_i,\, (v_i,t)\in\Gin\times (0,T).
\end{equation}
The vanishing of the other  two terms for each $v_i\in\Gin$, namely
\begin{equation}\label{kir2}
\sum_{j\in Inc_i}\nu_j\pd_j m(v_i,t)+ \pd_pH(v_i,\pd u)m_j(v_i,t)=0,
\end{equation}
gives the transition condition for $m$  at the vertices  $v_i\in\Gin$. Note that    \eqref{kir2} gives the conservation of the total flux of the density $m$ at the vertex $v_i$ (see   \cite{cg} for a similar condition). Summarizing, for $\nu:=\{\nu_j\}_{j\in J}$, we get the system
\begin{equation}\label{MFG}
    \left\{
    \begin{aligned}
    &-u_t-\nu \pd^2 u  +H(x, Du)=V[m]\qquad &(x,t)\in \G\times(0,T)\\[2pt]
    &m_t-\nu \pd^2 m -\pd(m\,\pd_p H(x, \pd u))=0 &(x,t)\in \G\times(0,T)\\[4pt]
    &\sum_{j\in Inc_i}\a_{ij}\nu_j\pd_j u(v_i,t)=0 &(v_i,t)\in\Gin\times (0,T)\\[2pt]
    &\sum_{j\in Inc_i}\nu_j\pd_j m(v_i,t)+ \pd_pH(v_i,\pd u)m_j(v_i,t)=0 &(v_i,t)\in\Gin\times (0,T)\\[4pt]
    &u_j(v_i,t)=u_k(v_i,t),\, \frac{m_j(v_i,t)}{\a_{ij}}= \frac{m_k(v_i,t)}{\a_{ik}} &j,k\in Inc_i,\, (v_i,t)\in\Gin\times (0,T)\\[4pt]
    &  u(x,T)=V_0[m(T)], \quad m(x,0)=m_0(x) &x\in\G
    \end{aligned}
    \right.
\end{equation}
with the normalization condition $\int_\G m(x)dx=1$.
The transition conditions (continuity and either Kirchhoff condition or conservation of total flux) for $u$ and $m$ give $d_{v_i}$ linear conditions  for each functions at a  $v_i\in\Gin$, hence  they univocally   determine the values $u_j(v_i,t)$ and $m_j(v_i,t)$, $j\in Inc_i$.

\vskip 8pt

\emph{Optimal control:}
We consider the \emph{planning problem} for a  MFG system,  i.e. we prescribe the initial and the terminal condition
for the distribution $m$ (see \cite{acc, ll}). Consider the functional
\begin{equation}\label{functional}
    \inf_{b,m} \int_0^T\int_\G\{L(x,b)m+W[m]\}dxdt
\end{equation}
subject to
\begin{equation}\label{constraint}
    \left\{
    \begin{aligned}
    &m_t-\nu \pd^2 m -\pd(bm)=0&(x,t)\in \G\times(0,T)\\
    &\sum_{j\in Inc_i}\nu_j \pd_jm(v_i,t)+b(v_i,t) m_j(v_i,t)=0&(v_i,t)\in \Gin\times (0,T)\\
    &\frac{m_j(v_i,t)}{\a_{ij}}=\frac{m_k(v_i,t)}{\a_{ik}}   &(v_i,t)\in\Gin\times (0,T)\,j,k\in Inc_i, \\
    & m(x,0)=m_0(x),\quad m(x,T)=m_T(x),& x\in\G.
    \end{aligned}
    \right.
\end{equation}
The problem of minimizing \eqref{functional}
under the constraints \eqref{constraint} is equivalent  to
\begin{equation}\label{functional_int}
    \inf_{b,m}\sup_{u} \int_0^T\int_\G\{L(x,b)m+W[m]-u(m_t-\nu \pd^2 m-\pd(mb))\} dxdt
\end{equation}
where $u$ is the multiplier.
We argue as in \cite[section 3.3]{acc} (see also \cite[section 2.5.1]{gs}); integrating by part in \eqref{functional_int}, taking into account the transition conditions in \eqref{constraint} and minimizing with respect to $b$, we obtain a minimum problem whose optimality conditions  give, at a formal level, a system similar to \eqref{MFG} with an initial-terminal condition for $m$ with $V=W'$.

\vskip 8pt

Similar considerations in both  the approaches can be used for deriving the stationary ({\it ergodic}) MFG system
\begin{equation}\label{MFGs}
    \left\{
    \begin{aligned}
     -&\nu \pd^2 u  +H(x, \pd u)+\rho  =V[m]\qquad &x\in \G\\
      &\nu \pd^2 m +\pd(m\,\pd_p H(x, \pd u))=0&x\in \G\\
        & \sum_{j\in Inc_i}\a_{ij}\nu_j\pd_j u(v_i)=0 & v_i\in\Gin\\
        &\sum_{j\in Inc_i}[ \nu_j\pd_j m(v_i)+\pd_pH_{j}(v_i,\pd_j u) m_j(v_i)]=0 & v_i\in\Gin\\
         &u_j(v_i)=u_k(v_i),\, \frac{m_j(v_i)}{\a_{ij}}= \frac{m_k(v_i)}{\a_{ik}} &j,k\in Inc_i,\, v_i\in\Gin\\
        &\int_\G u(x)dx=0,\quad\int_\G m(x)dx=1
    \end{aligned}
    \right.
\end{equation}
where $\rho\in\R$ is also an unknown.

In the rest of the paper we will only consider the stationary system \eqref{MFGs}. Moreover we will restrict
to the   case in which all the coefficients in  the transition condition for $u$ are   equal, i.e.
\begin{equation}\label{constant}
\a_{ij}=\a_{ik}\qquad\forall i\in I,\, j,k\in Inc_i.
\end{equation}
If \eqref{constant} is not satisfied, the function $m$ should be discontinuous at $v_i$ and this fact clearly involves additional difficulties. In fact it is well
known  that, in its standard definition, the domain of the Laplace operator on a network  is given by the $H^2(\G)$ functions (in particular, continuous) satisfying Kirchhoff condition at the vertices (\cite{fw,lu}). Moreover, the continuity condition at transition vertices seems to be a crucial ingredient for the comparison principle (see \cite{nvb,nic}).


\section{Main results}\label{sec3}

This section contains our main results on the solvability of HJB equations, FP equations and, above all, MFG systems on networks. To this end, we first introduce some assumptions.
Consider an Hamiltonian $H:\G\times  \R\to \R$, namely a collection of operators $(H_j)_{j\in J}$ with $H_j: [0,l_j]\times\R\to \R$. For some $\d$ and $C$ positive numbers, we assume
\begin{align}
    &H_j\in C^2([0,l_j]\times \R);\label{H1}\\
    &H_j(x,\cdot)\quad \text{is convex in $p$ for each $x\in [0,l_j]$};\label{H2}\\
    &\d |p|^2-C\le H_j(x,p)\le C|p|^2+C \quad\text{for  $(x,p)\in [0,l_j]\times \R$;}\label{H3}\\
    &\n=\{\n_j\}_{j\in J}, \quad \nu_j \in\R \quad\text{ with } 0<\nu_0:=\inf_{j\in J}\nu_j.\label{discount}
\end{align}
These assumptions will hold throughout this paper unless it is explicitly assumed in a different way. Let us note that no continuity condition for $H$ is required at the vertices and that, clearly, also the diffusion $\nu$ may present discontinuities at these points.

Let us now state our result for MFG systems, whose proof is contained in Section \ref{sec4}; in Sections \ref{HJBsubsect} and \ref{FPsubsect} we shall establish our result for HJB equations and respectively for FP equations.

\begin{Theorem}\label{existenceMFGs}
Assume \eqref{H1}-\eqref{discount} and that $V$ is a local $C^1$ coupling, namely
\begin{equation}\label{V}
V[m] (x)= V(m(x)) \textrm{ with }V\in C^1([0,+\infty)).
\end{equation}
Then, there exists a solution $(u,m,\rho)\in C^2(\G)\times C^2(\G)\times\R$ to
\begin{equation}\label{MFGzero}
    \left\{
    \begin{aligned}
     -&\nu \pd^2 u  +H(x, \pd u)+\rho =V[m]\qquad &x\in \G\\
    &\nu \pd^2 m +\pd(m\,\pd_p H(x, \pd u))=0&x\in \G\\
    & \sum_{j\in Inc_i}\nu_j\pd_j u(v_i)=0  &v_i\in\Gin\\
    &\sum_{j\in Inc_i}[\nu_j\pd_j m(v_i)+ \pd_pH_{j}(v_i,\pd_j u) m_j(v_i)]=0 &  v_i\in\Gin\\
    &u_j(v_i)=u_k(v_i),\,  m_j(v_i)=  m_k(v_i) &j,k\in Inc_i,\, v_i\in\Gin\\
    &\int_\G u(x)dx=0,\quad \int_\G m(x)dx=1,\qquad m\ge 0.
    \end{aligned}
    \right.
\end{equation}
Moreover if
\begin{equation}\label{V1}
    \int_\G (V[m_1]-V[m_2])(m_1-m_2)dx\le 0  \Rightarrow \text{ $m_1=m_2$,}
\end{equation}
then the solution is unique.
\end{Theorem}
\begin{Remark}
As already pointed out in the introduction, this result can be easily adapted to the case of networks having a boundary by imposing Neumann or Dirichlet boundary condition.
\end{Remark}

\subsection{On the  Hamilton-Jacobi-Bellman equation}\label{HJBsubsect}

This section is devoted to the ergodic HJB problem: find $(u,\rho)\in C^2(\G)\times \R$ such that
\begin{equation}\label{hjberg}
\left\{\begin{array}{lll}
    & -\nu \pd^2 u  +H(x, \pd u)+\rho =f(x),\qquad &x\in \G\\[3pt]
    &u_j(v_i)=u_k(v_i),\quad \ds\sum_{j\in Inc_i}\nu_j\pd_j u(v_i)=0&j,k\in Inc_i,\quad v_i\in\Gin
    \end{array}
    \right.
\end{equation}
with the normalization condition
\begin{equation}\label{hjbnorm}
    \int_\G u(x)dx=0.
\end{equation}

\begin{Theorem}\label{existhjb}
 Assume \eqref{H1}-\eqref{discount} and $f\in C^{0,\a}(\G)$ for some $\a\in(0,1)$. Then, there exists a  unique couple $(u,\rho)\in C^2(\G)\times \R$ satisfying   \eqref{hjberg}-\eqref{hjbnorm}. Moreover $u\in C^{2,\alpha}(\G)$ and there holds
 \begin{equation}\label{esthjb}
    \|u\|_{C^{2,\a}(\G)}\le C,\qquad |\rho|\leq \max_\G|H(\cdot,0)-f(\cdot)|
 \end{equation}
with $C$ depending only  on $\|f\|_{C^{0,\a}}$ and the constants in \eqref{H3}-\eqref{discount}.
\end{Theorem}
\begin{Proof}
Proposition \ref{bmp} ensures that, for any $\l\in(0,1)$, there exists a solution $u_\l\in C^{2,\a}(\G)$ to
\begin{equation}\label{hjb}
\left\{\begin{array}{lll}
    & -\nu \pd^2 u  +H(x, \pd u)+\l u =f(x),\qquad &x\in \G\\[3pt]
    &u_j(v_i)=u_k(v_i),\quad \sum_{j\in Inc_i}\nu_j\pd_j u(v_i)=0&j,k\in Inc_i,\quad v_i\in\Gin.
    \end{array}
    \right.
\end{equation}
We want to pass to the limit for $\l\to 0$ in \eqref{hjb}. We first observe that if $C_0$ is a constant such that $\max_\G|H(\cdot,0)-f(\cdot)|\le C_0$, then the functions $\underline u$, $\overline u$ defined by $\underline u(x)=-C_0/\l$,  $\overline u(x)=C_0/\l$ for any $x\in \G$, are respectively a sub- and a supersolution of \eqref{hjb}.
By Proposition \ref{maxprin} we get
\begin{equation}\label{limit0}
 -C_0\le \l u_\l(x)\le C_0\quad \text{for any $x\in\G$.}
\end{equation}
Now, let us introduce the function $w_\l:=u_\l-\min_\G u_\l$; it is a $C^{2,\a}$ solution to
\begin{equation}\label{3.14}
-\nu \pd^2 w_\l  +H(x,\pd w_\l)+\l u_\l  =f(x)\qquad x\in \G
\end{equation}
with the same continuity and Kirchhoff transition conditions as in \eqref{hjb}. We claim
\begin{equation}\label{cl:goodest}
\|\pd w_\l\|_{L^2(\G)}\leq C_1
\end{equation}
for some constant $C_1$ depending only on $\|f\|_{C^{0,\a}}$ and the constants in \eqref{H3}-\eqref{discount} (in particular, independent of $\l$). Indeed, integrating equation \eqref{3.14} on $\G$ (i.e. using $\phi=1$ as test function for $w_\l$), we get
\[
\int_\G H(x, \pd w_\l)\, dx +\int_\G(\l u_\l)\, dx = \int_\G f \, dx.
\]
By assumption \eqref{H3} and estimate \eqref{limit0}, we infer
\[
\d \int_\G |\pd w_\l|^2\, dx\leq \int_\G(C-\l u_\l+f)\, dx\leq (C+C_0+\|f\|_\infty)|\G| \]
which amounts to our claim \eqref{cl:goodest}.
We claim now that
\begin{equation}\label{cl:goodest1}
\| w_\l\|_{C^{2,\a}(\G)}\leq  C_2
\end{equation}
for some constant $C_2$ with the same feature of $C_1$. To this end, we note that, since $w_\l$ is a classical solution to \eqref{3.14}, by \eqref{H3}-\eqref{discount} and \eqref{limit0}, there holds
\[
\nu_0|\pd^2 w_\l|\leq |H(x, \pd w_\l)|+|\l u_\l|\leq C(|\pd w_\l|^2 +1)+C_0
\]
and, by \eqref{cl:goodest}
\begin{equation}\label{cl:goodest2}
\|\pd^2 w_\l\|_{L^1(\G)}\leq C_3
\end{equation}
for some constant $C_3$ sharing the same features of $C_1$.
Taking into account \eqref{cl:goodest} and \eqref{cl:goodest2}, (possibly increasing $C_3$) we infer $\|\pd w_\l\|_{C^{0,\a}(\G)}\leq C_3$; using again \eqref{3.14} we accomplish the proof of our claim \eqref{cl:goodest1}.

Possibly passing to a subsequence, we may assume that, as $\l\to 0^+$, the sequence $\{w_\l\}_\l$ converges to some function $u\in C^{2,\a}(\G)$ (observe that $u$ still verifies the continuity and the Kirchhoff conditions) and that $\{\l \min_\G u_\l\}_\l$ converges to some constant $\rho$.
Passing to the limit in \eqref{3.14}, we get that the couple  $(u,\rho)$ satisfies
 \eqref{hjberg}. Possibly adding a constant to $u$, we also get \eqref{hjbnorm}.\\
 To show the uniqueness of $\rho$, assume that there exist two solutions $(u_i,\rho_i)$, $i=1,2$, of \eqref{hjberg} and let $x_0$ be a maximum point of $u_1-u_2$. If  $x_0\in e_j$, we have $\pd_j u_1(x_0)=\pd_j u_2(x_0)$ and $\pd^2_j u_1(x_0)\le \pd^2_j u_2(x_0)$. Hence using the equation  we conclude that $\rho_2\le \rho_1$. If $x_0=v_i$, there holds $\pd_j u_1(x_0)\leq \pd_j u_2(x_0)$  for any $j\in Inc_i$. In fact, we have: $\pd_j u_1(x_0)=\pd_j u_2(x_0)$  for any $j\in Inc_i$; indeed, assuming by contradiction $\pd_k u_1(x_0)< \pd_j u_2(x_0)$  for some $k\in Inc_i$, we get $\sum_{j\in Inc_i}\nu_j \pd_j u_1(x_0)<\sum_{j\in Inc_i}\nu_j \pd_j u_2(x_0)$ which contradicts the Kirchhoff condition. Hence
\[
\nu_j\pd^2_j(u_2-u_1)(x_0)\geq H(x_0, \pd_j u_2(x_0))-H(x_0, \pd_j u_1(x_0))+\rho_2-\rho_1=\rho_2-\rho_1
\]
which, together with  $\pd_ju_1(x_0)=\pd_j u_2(x_0)$, gives a contradiction for $\rho_2>\rho_1$. Hence  $\rho_2\le \rho_1$ and by symmetry $\rho_2=\rho_1$.\\
Having proved the uniqueness of $\rho$,   the  uniqueness of a solution to \eqref{hjberg}-\eqref{hjbnorm} can be proved as in    \cite[Corollary 3.1]{cms}. Finally, since there exists a unique solution to \eqref{hjberg}  we conclude  that all the sequence $(w_\l,\l u_\l)$ converges to $(u,\rho)$.
\end{Proof}
\subsection{On the Fokker-Planck equation}\label{FPsubsect}
This section is devoted to the following problem
\begin{equation}\label{fpCLA}
\left\{
\begin{array}{ll}
\nu \pd^2 m +\pd(b(x)\,m)=f(x)\qquad & x\in\G\\[4pt]
 m_j(v_i)=  m_k(v_i),\quad \ds\sum_{j\in Inc_i} [ b(v_i) m_j(v_i)+\nu_j\pd_j m(v_i)]=0 & j,k\in Inc_i,\, v_i\in\Gin
      \end{array}
      \right.
\end{equation}
 with the normalization conditions
\begin{equation}\label{normalization}
m\ge 0, \qquad \int_\G m(x)dx=1.
\end{equation}

\begin{Definition}
\begin{itemize}
  \item[(i)] A  strong  solution of \eqref{fpCLA} is a function $m\in C^2(\G)$ which satisfies \eqref{fpCLA} in pointwise sense.
  \item[(ii)]  A weak solution of \eqref{fpCLA} is a function $m\in H^1(\G)$ such that
\begin{equation}\label{weakCLA}
 \sum_{j\in J}\int_{e_j} (\nu_j\pd_j m+b(x)m)\pd_j \phi\, dx+\int_\G f \phi dx=0\qquad\forall\phi\in H^1(\G).
\end{equation}
\end{itemize}
\end{Definition}
\begin{Remark} By standard arguments, one can easy check that if $m\in C^2(\G)$ is a weak solution to \eqref{fpCLA}, then it is also a strong solution to \eqref{fpCLA}.
\end{Remark}
\begin{Theorem}\label{existfpCLA}
Assume $b\in C^1(\G)$. Then, there exist a unique weak  solution $m$ to  \eqref{fpCLA}-\eqref{normalization} with $f=0$ and it verifies
 \begin{equation}\label{estfpCLA}
    \|m \|_{H^1}\le C,\qquad 0< m(x)\leq C
\end{equation}
where the constant $C$ depends only on $\|b\|_{\infty}$ and $\nu$.
Moreover $m\in C^2(\G)$ and it is also a strong solution to \eqref{fpCLA}.
 \end{Theorem}
\begin{Proof}
We shall proceed adapting the arguments of \cite[Theorem II.4.3]{ben}. By Proposition~\ref{ben4.2}, there exists a unique (up to multiplicative constant) solution to problem~\eqref{fpCLA}. So we only have to prove that this family of solutions contains a (unique) function~$m$ satisfying~\eqref{normalization} which moreover will verify~\eqref{estfpCLA}.\\
For $\l\in(1,+\infty)$ and for any $\phi\in L^\infty(\G)$, we set
\[
U_\l(t,x):=U(\l t,x)\qquad \textrm{for }(t,x)\in(0,1)\times\G
\]
where $U$ is the solution of the parabolic Cauchy problem \eqref{dualpara} with $\psi=\phi$ (see Lemma~\ref{lemma:pdfp}). Observe that $U_\l$ solves
\begin{equation}\label{ben425}
\pd_t U_\l- \l \nu \pd^2 U +\l b \pd U=0 \qquad \textrm{in }(0,1)\times\G
\end{equation}
with the same transition conditions and initial datum of~\eqref{dualpara}. We claim that
\begin{equation}\label{ben426}
\textrm{$U_\l$ are uniformly bounded in $L^2(0,1;H^1(\G))\cap L^\infty((0,1)\times\G)$.}
\end{equation}
Indeed, since $\pm\|\phi\|_\infty$ are respectively a super- and a subsolution to~\eqref{ben425}, we get
\begin{equation}\label{ben426bis}
|U_\l(t,x)|\leq \|\phi\|_\infty\qquad \textrm{a.e. in }(0,1)\times\G.
\end{equation}
In other words, we get that the functions $U_\l$ are uniformly bounded in $L^\infty((0,1)\times\G)$ and, in particular, in $L^2((0,1)\times\G)$.
On the other hand, using $U_\l$ as test function for problem~\eqref{ben425}, we get
\begin{eqnarray*}
\int_\G U^2_\l(\t,x) \, dx&+&\l \sum_j\nu_j\iint_{(0,1)\times e_j}(\pd_j U_\l)^2\, dxdt\\
&=& \int_\G \phi^2(x) \, dx -\l \iint_{(0,1)\times  \G} b\pd U_\l U_\l\, dxdt\\
&\leq& \int_\G \phi^2(x) \, dx +\l\|b\|_\infty \iint_{(0,1)\times \G} |\pd U_\l||U_\l|\, dxdt.
\end{eqnarray*}
Applying Cauchy inequality to the last term (recall that $\G$ has finite measure), by \eqref{ben426bis}, we get
\[
\l \sum_j\nu_j\iint_{(0,1)\times e_j}(\pd_j U_{\l})^2\, dxdt\leq c(\l+1)
\]
for some constant $c$ independent of $\l$. We infer that $\pd U_\l$ are uniformly bounded in $L^2((0,1)\times\G)$; thus our claim~\eqref{ben426} is completely proved.

The property~\eqref{ben426} yields that there exists a subsequence of $\{U_\l\}$ (that we still denote by $U_\l$) such that
\[
U_\l \to \xi \quad \textrm{in $L^\infty((0,1)\times\G)$ weak-$*$ and in $L^2(0,1;H^1(\G))$ weak as $\l\to +\infty$}.
\]
Let us now use $\b\th$ as test function for~\eqref{ben425}, with $\b\in C^\infty_0((0,1))$ and $\th\in C^\infty(\G)$ (recall: this means that $\th\in C^0(\G)$ and $\th\in C^\infty(e_j)$ for every $j\in J$); we obtain
\[
\frac 1\l \iint_{[0,1]\times \G} U_\l\b'\th dx\,dt+\int_{[0,1]}\b\sum_j\int_{e_j}\left(\nu_j\pd_j U_\l \pd_j \th + b \pd_j U_\l \th \right)dx\,dt =0.
\]
Passing to the limit as $\l\to+\infty$, we get
\[
\int_{[0,1]}\b\sum_j \int_{e_j}\left(\nu_j\pd_j \xi \pd_j \th+
b \pd_j \xi \th\right)dx\,dt=0.
\]
By the arbitrariness of $\b$ we get that, for a.e. $t\in(0,1)$, there holds
\[
\sum_j\int_{e_j}\left(\nu_j\pd_j \xi \pd_j \th+b \pd_j \xi \th\right)dx=0;
\]
namely, $\xi$ is a weak solution to
\begin{equation*}
\left\{
\begin{array}{lll}
&-\nu \pd^2 \xi +b(x)\pd \xi=0\qquad & x\in\G\\[4pt]
&\xi_j(v_i)=\xi_k(v_i),\quad\ds\sum_{j\in Inc_i} \nu_j\pd_j \xi(v_i)=0 &j,k\in Inc_i,\, v_i\in\Gin.
      \end{array}
      \right.
\end{equation*}
The maximum principle for this problem (see~\cite[Theorem 2.1]{nic}) ensures that the function $\xi$ is independent of $x$, namely $\xi=\xi(t)$.

On the other hand, using the function~$m$ introduced in Proposition~\ref{ben4.2} as test function for~\eqref{ben425}, we infer
\[
\int_\G U_\l(\t,x)m(x)\, dx = \int_\G \phi(x)m(x)\, dx;
\]
as $\l\to+\infty$, we get
\[
\xi(t)\int_\G m(x)\, dx = \int_\G \phi(x)m(x)\, dx.
\]
By the arbitrariness of $\phi$ (recall also that $m$ cannot be identically zero because it belongs to a $1$-dimensional family), we deduce that $\int mdx$ cannot be zero. Moreover, we also infer that $\xi$ is independent of $t$, namely $\xi$ is a constant. By Proposition~\ref{ben4.2}, we can choose $m$ such that $\int mdx=1$. In conclusion, the  last equality reads as
\[
\xi = \int_\G \phi(x)m(x)\, dx
\]
for every $\phi\in L^\infty(\G)$ (clearly, $\xi$ depends on $\phi$).
By Lemma \ref{lemma:pdfp}, a standard application of the Kantorovich-Vulikh theorem (see, \cite[Theorem 6.3]{mug} and the subsequent discussion) ensures that the semigroup associated to the Cauchy problem~\eqref{dualpara} has a {stricly positive} integral kernel. Therefore, we may accomplish the proof following the same arguments of~\cite{ben} and of \cite[Lemma 2.3]{bf}.\\
Finally, let us prove  that $m$ belongs to $C^2(\G)$; fix $j\in J$ and rewrite the equation as
\begin{equation}\label{eq2}
    \nu \pd^2_j m=-m\pd_j b-b\pd_j m
\end{equation}
Since $b\in C^1(\G)$ and $m\in H^1(\G)$, it follows  that $m\in H^2(\G)$, hence $\pd_j m$ is continuous.  Therefore by \eqref{eq2}
we conclude that $m\in C^2(\G)$ and it is also  a strong solution to \eqref{fpCLA}.
\end{Proof}

\subsection{Proof of Theorem \ref{existenceMFGs}}\label{sec4}

\begin{Proofc}{Proof of Theorem \ref{existenceMFGs}}
Consider the set  $\cK=\{\mu\in C^{0,\a}(\G):\, \int_\G \mu dx=1\}$ and observe that $\cK$ is a closed subset of the Banach space $C^{0,\a}(\G)$. We define an operator $T:\cK\to\cK$ according to the scheme
\begin{equation}\label{mapT}
\mu \to u \to m
\end{equation}
as follows. Given $\mu\in \cK$, we  solve the problem \eqref{hjberg}-\eqref{hjbnorm} with $f(x)=V[\mu](x)$ for the unknowns $u$ and $\rho$, which are uniquely defined by Theorem \ref{existhjb}. Then, for $u$ given, we seek a function $m$ which solves problem \eqref{fpCLA}-\eqref{normalization} with $b(x)=\pd_p H(x,\pd u)$. By  Theorem \ref{existfpCLA}, the function $m$ is univocally defined and we set $m=T( \mu)$. We claim that
\begin{equation}\label{claim1}
\textrm{the map $T$ is continuous with  compact image.}
\end{equation}
 Let $\mu_n,\mu \in \cK$ be such that $\|\m_n-\m\|_{C^{0,\a}}\to 0$ for $n\to \infty$ and let $(u_n,\rho_n)$, $(u,\rho)$ be the solutions of \eqref{hjberg}-\eqref{hjbnorm} corresponding to $f(\cdot)=V(\mu_n(\cdot))$ and, respectively, $f(\cdot)=V(\mu(\cdot))$.  By the estimate \eqref{esthjb}, (possibly passing to a subsequence) $u_n$ converges in $C^2(\G)$ to a function $\bar u$ and $\rho_n$ converges to some constant $\bar \rho$. Since $V[\m_n]\to V[\m]$ in $C^{0,\a}$ and since the transition and  the boundary conditions pass to the limit by the $C^2$-convergence we get that $\bar u$ is a solution of \eqref{hjberg}-\eqref{hjbnorm} with $f(x)$ and $\rho$ replaced respectively by $V[m](x)$ and $\bar\rho$. By the uniqueness of the  solution to \eqref{hjberg}-\eqref{hjbnorm} we get that $\bar u=u$ and $\rho=\bar \rho$; moreover,  all the sequence $\{(u_n,\rho_n)\}$ converges to $(u,\rho)$. Let $m_n$ and $m$ be respectively the solutions of \eqref{fpCLA}-\eqref{normalization} with $f=0$ for $b=\pd_p H(x,\pd u_n)$ and for $b=\pd_p H(x,\pd u)$.
By the estimate \eqref{estfpCLA}, the functions $m_n$ are equibounded in $H^1(\G)$. Since $\pd_p H(x,\pd u_n)$ uniformly converges to $\pd_p H(x,\pd u)$, passing to the limit in the weak formulation yields that (possibly passing to a subsequence) $m_n$ converges to a solution $\bar m$ to \eqref{fpCLA}-\eqref{normalization} with $b(x)=\pd_p H(x,\pd u)$. Theorem \ref{existfpCLA} entails: $m=\bar m$; hence, the whole sequence $\{m_n\}$ converges to $m$.\par
To prove the compactness of the image of $T$, consider  a sequence  $\m_n$ such that $\|\m_n\|_{C^{0,\a}}\le 1$ and let $u_n$ and $m_n$ the functions obtained according to the scheme \eqref{mapT}. By \eqref{esthjb}, $\|u_n\|_{C^{2,\a}}$ is uniformly bounded and therefore by \eqref{estfpCLA} also $\|m_n\|_{{C^{0,\a}}}$ is uniformly bounded. As in the proof of Proposition \ref{existfpCLA}, we get an uniform bound on the $H^2$-norm of $m_n$ for any $n$. By the compact immersion of $H^2(\G)$ in $H^1(\G)$, we get that the sequence $m_n$ is
compact in $H^1$ and therefore in $C^{0,\a}$. Hence, claim \eqref{claim1} is completely proved.\par
We can therefore conclude  by the Leray-Schauder fixed point Theorem that the map
 $T$ admits a fixed point , i.e. a solution of system \eqref{MFGzero}. Moreover this solution is also smooth by the regularity results in Theorems \ref{existhjb} and \ref{existfpCLA}.

Finally the uniqueness of the solution to \eqref{MFGzero} under the assumption \eqref{V1} follows by a standard argument in MFG theory adapted to the networks (see \cite{ll}).  We assume that there exists two solutions $(u_1,m_1,\rho_1)$ and $(u_2,m_2,\rho_2)$ of \eqref{MFGzero}. We set $\bar u=u_1-u_2$, $\bar m=m_1-m_2$,
$\bar \rho=\rho_1-\rho_2$ and we write the equations for $\bar u$, $\bar m$
\[\left\{\begin{array}{ll}
    -  \nu\pd^2 \bar u+H(x, \pd  u_1)-H(x, \pd  u_2)+\bar \rho-(V[m_1] -V[m_2])=0 \\[2pt]
       \nu\pd^2 \bar m +\pd (m_1 \pd_p H (x,\pd u_1)-m_2\pd_p H (x,\pd u_2))=0 \\[2pt]
      \sum_{j\in Inc_i}\nu_j\pd_j \bar u(v_i)=0\\[2pt]
      \sum_{j\in Inc_i}\nu_j\pd_j \bar m(v_i)+(m_1 \pd_p H(v_i,\pd u_1)-m_2 \pd_p H(v_i,\pd u_2))=0\\[2pt]
     \int_\G \bar m dx=0,\quad \int_\G \bar u dx=0
\end{array}\right.\]
Multiplying the equation for $\bar m$ by $\bar u$ and integrating over $e_j$, we get
\begin{equation}\label{U4}
\begin{split}
 &\int_{e_j} \left[-\nu_j \pd_j \bar u\pd_j\bar m - \big( m_1\pd_p H_j(x,\pd u_1)-m_2\pd_p H_j(x,\pd u_2)\big) \pd_j \bar u(x)\right] dx \\
 &+\Big [ \bar u_j(\nu_j\pd_j \bar m + m_1\pd_p H_j(x,\pd u_1)-m_2\pd_p H_j(x,\pd u_2))\Big]_{0}^{ l_j} =0.
\end{split}
\end{equation}
Multiplying the equation for $\bar u$  by $\bar m$ and integrating over $e_j $, we get
\begin{equation}\label{U5}
\begin{split}
 &\int_{e_j}\nu_j \pd_j \bar u\pd_j\bar m +\big[H_j(x, \pd_j  u_1)-H_j(x, \pd_j  u_2)+\bar \rho -(V[m_1] -V[m_2])\big]\bar m_j dx\\
 &+\Big [\nu_j\bar m \pd_j \bar u\Big]_{0}^{ l_j}=0
\end{split}
\end{equation}
Adding \eqref{U4} to \eqref{U5}, summing over $j\in J$, regrouping the terms corresponding to a same vertex $v_i$ and  taking into account the transition and the normalization conditions for $\bar u$ and $\bar m$
we get
\begin{align*}
    &\sum_{j\in J} \int_{e_j}(m_1-m_2) (V[m_1] -V[m_2])dx +\\
    & \sum_{j\in J} \int_{e_j} m_1\big[H_j(x,\pd_j u_2) - H_j(x,\pd_ju_1) -\pd_p H_j(x,\pd_j u_2)\pd_j( u_2-u_1)\big]dx +\\
     &\sum_{j\in J} \int_{e_j} m_2\big[(H_j(x,\pd_j u_1) - H_j(x,\pd_ju_2) -\pd_p H_j(x,\pd_j u_1)\pd_j( u_1-u_2)\big]dx=0.
    \end{align*}
Since each of the three terms in the previous identity is non-negative, it follows that it
must vanish. By \eqref{V1} we get $m_1=m_2$. By the uniqueness of the solution to \eqref{hjberg} we finally get $u_1=u_2$  and $\rho_1=\rho_2$.
\end{Proofc}


\appendix
\section{Appendix}\label{appendix}
\subsection{Auxiliary  results for  HJB equation}\label{appA1}
In this section, we study the following semilinear problem
\begin{equation}\label{bpm1}
\left\{\begin{array}{lll}
    & -\n\pd^2 u+H(x,u,\pd u)+\l  u=0,\qquad &x\in \G\\[4pt]
    &u_j(v_i)=u_k(v_i),\quad\ds\sum_{j\in Inc_i}\nu_j\pd_j u(v_i)=0&j,k\in Inc_i,\, v_i\in\Gin.
    \end{array}
    \right.
\end{equation}
As far as we know, this problem has not been tackled before; we shall establish existence, regularity and uniqueness (via comparison principle).
\begin{Definition}\label{weaksolhjb}\hfill
\begin{itemize}
\item[(i)] A  strong  solution of \eqref{bpm1} is a function $u\in C^2(\G)$ which satisfies \eqref{bpm1} in pointwise sense.
\item[(ii)]  A weak solution of \eqref{bpm1} is a function $u\in H^1(\G)$ such that
\begin{equation}\label{weak}
 \sum_{j\in J}\int_{e_j} (\nu_j\pd_j u\pd_j\phi+H(x,u,\pd_j u)\phi+\l u\phi)dx=0\qquad\text{for any  $\phi\in H^1(\G)$.}
\end{equation}
\end{itemize}
\end{Definition}
\begin{Remark}
One can easily check that, if $u\in C^2(\G)$ is a weak solution of \eqref{bpm1}, then it is also a strong solution.
\end{Remark}
Let us now state our existence result
\begin{Proposition}\label{bmp}
Assume
\begin{align}
&
\begin{array}{l}
   \text{$H_j(\cdot,r,p)$ is measurable in $x$, for any $(r,p)\in \R\times\R$ and}\\
   \text{$H_j(x,\cdot, \cdot)$ is continuous  in $(r,p)$, for a.e. $x\in (0,l_j)$}
   \end{array}\label{regolarita}\\[4pt]
   & |H(x,r,p)|\le C_0+b(|r|)|p|^2\qquad \text{for a.e. $(x,r,p)\in\G\times\R\times\R$}\label{subquad}\\
   & H(x,r,p)\text{ is not decreasing in $r$ for a.e. $(x,r,p)\in\G\times\R\times\R$}\label{monoton}\\
   &\l>0,\, \nu_j \in\R \quad\text{ with } 0<\nu_0:=\inf_{j\in J}\nu_j\label{zeroorder}
\end{align}
where $C_0>0$, $ b:\R\to\R$ is an increasing function.
 Then there exists a  weak solution to  \eqref{bpm1}. Moreover
   \[ \|u\|_{H^1}\le C\]
with $C$ depending on $C_0$, $\l$ and   $\n_0$.

Moreover, if $H$ belongs to $C^{0,\a}(\G\times \R \times \R)$ for some $\a\in(0,1)$, then solution $u$ belongs to $C^{2,\a}(\G)$ with
\begin{equation}\label{strongest}
\|u\|_{C^{2,\a}}\leq C_1(1+\|u\|_{H^1})
\end{equation}
where $C_1$ is a constant depending only on $C_0$, $b$ and $\nu_0$.
\end{Proposition}
\begin{Lemma}\label{bpmL1}
Assume \eqref{regolarita}, \eqref{monoton} and, for some $C_H>0$,
\[
|H(x,r,p)|\le C_H \quad\text{for  $(x,r,p)\in \G\times\R\times\R$}.
    \]
Then there exists a weak solution $u$ to \eqref{bpm1}. Moreover
\begin{equation}\label{weakmax}
  \|u\|_\infty\le\frac{C_H}{\l}.
\end{equation}
\end{Lemma}
\begin{Proof}
 Define  a  map $T:H^1 (\G)\to H^1(\G)$ by taking the weak solution  $u=T(v)$ of
\begin{equation}\label{E1}
\left\{
\begin{array}{lll}
&- \nu \pd^2 u+ \l u=-H(x,v,\pd v)\qquad & x\in\G\\[4pt]
&u_j(v_i)=u_k(v_i),\quad\ds\sum_{j\in Inc_i}\nu_j\pd_j u(v_i)=0&j,k\in Inc_i,\, v_i\in\Gin.
      \end{array}
      \right.
\end{equation}
(note that existence of a weak solution to \eqref{E1}  follows by  the theory of sesqui-linear forms, see for instance \cite{nic}). Standard estimates implies that $T$ is continuous with  compact image,  hence  by the Schauder's Theorem  it admits a fixed point which is a   weak solution to \eqref{bpm1}.

Even though the proof of estimate \eqref{weakmax} is standard, for completeness we sketch the argument. Let $G:\R\to\R$ be a smooth function such that $G(t)=0$ for $t\in (-\infty,0]$ and $G$ strictly increasing   for $t\in (0,\infty)$.
Set $K=C_H/\l$ and   $\phi=G(u-K)$. Then $\phi\in H^1(\G)$ and by taking $\phi$ as test function in \eqref{weak} we get
\[
  0=\sum_{j\in J}\int_{e_j} [\nu_j(\pd_j u)^2G'(u-K)+(H_j(x, u,\pd_j u )+\l K)G(u-K)+\l(u-K)G(u-K)]dx
\]
 Since $ H_j(x,u,\pd u)+\l K\ge 0$ and $G(u-K)\ge 0$ a.e. on $\G$, then
 \[\sum_{j\in J}\int_{e_j} \l(u-K)G(u-K)dx\le 0.\]
and by $tG(t)\ge 0$ for $t\in\R$, it follows that $(u-K)G(u-K)=0$ a.e. on $\G$, hence $u\le K$ a.e. in $\G$.
\end{Proof}
\begin{Proofc}{Proof of Proposition \ref{bmp}}
The proof of the existence result follows exactly the same argument of the proof of \cite[Theorem 2.1]{bmp81} replacing their steps $1$ and $2$ with Lemma \ref{bpmL1}. In fact, in the weak formulation of \eqref{hjb}, given in Definition \ref{weaksolhjb}, the transition condition is transparent and  the estimates necessary to prove the result are obtained using the weak formulation \eqref{weak} which is the same of the problem posed in an Euclidean domain.

Consider now $H\in C^{0,\a}(\G\times \R \times \R)$. We already know that $u\in C^0(\G)$  and we have only to show that $u_{j}\in C^{2,\a}(0,l_j)$  for any $j\in J$ (recall that $u_{j}$ is the restriction of $u$ to the edge $e_j$).  For $j$ fixed, the equation for $u$ (in distributional sense) is
\begin{equation}\label{eq1}
-\nu_j\pd^2_j u=-\l u_{j}-H(x,u_j,\pd_j u)\quad \textrm{in }(0,l_j).
\end{equation}
Since $u_{j}\in C^0([0,l_j])$ and, by \eqref{subquad}, $H(\cdot,u_j(\cdot),\pd_j u(\cdot))\in L^1((0,l_j))$, by \eqref{eq1} we immediately get $u_j\in W^{2,1}((0,l_j))$
and therefore $\pd_j u\in L^p((0,l_j))$, for any $p\ge 1$. We deduce  $H(\cdot,u_j(\cdot),\pd_j u(\cdot))\in L^p((0,l_j))$ and in particular $u_{j}\in W^{2,p}((0,l_j))$.
Hence  $\pd_j u\in C^{0,\a}((0,l_j))$  and again by \eqref{eq1} we get the statement. Moreover, the estimate \eqref{strongest} easily follows from \eqref{eq1}.
\end{Proofc}
\begin{Proposition} \label{maxprin}
Assume \eqref{regolarita} and \eqref{monoton}-\eqref{zeroorder}. Let the functions $u_1, u_2\in C^2(\G)$ satisfy
\begin{equation}\label{bmpeq10}
\left\{\begin{array}{ll}
-\n\pd^2 u_1+H(x,u_1,\pd u_1)+\l u_1\ge -\n\pd^2 u_2+H(x,u_2,\pd u_2)+\l u_2 &x\in \G\\[4pt]
 \sum_{j\in Inc_i}\nu_j\pd_j u_1(v_i)\le \sum_{j\in Inc_i}\nu_j\pd_j u_2(v_i) &v_i\in\Gin.
\end{array}\right.
\end{equation}
Then, $u_1\ge u_2$ on $\G$.
\end{Proposition}
\begin{Proof}
We argue by contradiction assuming $\max_{\G}(u_2-u_1)=:\d>0$.
Let $x_0$ be a point where $u_2-u_1$ attains its maximum. The point~$x_0$ either belongs to some edge or it is a vertex.
Assume that $x_0$ belongs to some edge $e_j$. By their regularity, the functions $u_1$ and $u_2$ fulfill
\[
u_2(x_0)=u_1(x_0)+\d,\qquad \pd_j u_2(x_0)=\pd_j u_1(x_0),\qquad \pd_j^2 u_2(x_0)\leq \pd_j^2 u_1(x_0).\]
In particular, we deduce
\begin{align*}
&-\n_j\pd_j^2u_1(x_0)+H(x_0,u_1(x_0),\pd_j u_1(x_0))+\l  u_1(x_0)\\
&\qquad\qquad \leq -\n_j\pd_j^2u_2(x_0)+H(x_0,u_2(x_0),\pd_j u_2(x_0))+  \l(u_2(x_0)-\d)\\
&\qquad\qquad<
-\n_j\pd_j^2u_2(x_0)+H(x_0,u_2(x_0), \pd_j u_2(x_0))+\l u_2(x_0)
\end{align*}
which contradicts the first relation in~\eqref{bmpeq10}.
Assume that $x_0=v_i$ for some $v_i\in\Gin$. Being regular, the functions $u_1$ and $u_2$ fulfill $\pd_j u_2(v_i)\leq \pd_j u_1(v_i)$. We claim $\pd_j u_2(v_i)= \pd_j u_1(v_i)$ for each $j\in Inc_i$. In order to prove this equality we proceed by contradiction and we assume that $\pd_j u_2(v_i)< \pd_j u_1(v_i)$ for some $j\in Inc_i$. In this case we get $\sum_{j\in Inc_i}\nu_j\pd_j u_2(v_i)< \sum_{j\in Inc_i}\nu_j\pd_j u_1(v_i)$ which contradicts the second inequality in \eqref{bmpeq10}; therefore, our claim is proved.
Moreover, since $u_1(v_i)=u_2(v_i)-\d$, we deduce
\begin{eqnarray*}
H(v_i,u_1(v_i),\pd_j u_1(v_i))+\l u_1(v_i)&=&H(v_i,u_2(v_i)-\d, \pd_j u_2(v_i))+\l(u_2(v_i)-\d)\\
&<&H(v_i, u_2(v_i),\pd_j u_2(v_i))+\l u_2(v_i).
\end{eqnarray*}
Taking into account the regularity of $H$ and of $u_i$ ($i=1,2$), we infer that in a sufficiently small neighborhood $B_\eta(v_i)$ there holds
\[H(x,u_1(x), \pd u_1(x))+\l u_1(x) <H(x,u_2(x),\pd u_2(x))+\l u_2(x).\]
This inequality and the first relation in \eqref{bmpeq10} entail
\begin{multline*}
\n_j\pd^2_j(u_2-u_1)(v_i) \geq H(v_i, u_2(v_i),\pd_j u_2(v_i))-H(x,u_1(v_i), \pd_j u_1(v_i))\\+\l(u_2(v_i)-u_1(v_i))>0
\end{multline*}
which, together with $\pd_j u_2(v_i)= \pd_j u_1(v_i)$, contradicts that $u_2-u_1$ attains a maximum in $x_0=v_i$.
\end{Proof}

\subsection{Auxiliary  results for  the Fokker-Planck equation}
\begin{Proposition}\label{ben4.2}
Assume that  $f\in C^0(\G)$ and $b\in C^1(\G)$.  The problem~\eqref{fpCLA} with $f=0$ admits a unique (up to multiplicative constant) weak solution.
Moreover, for $f\ne 0$, the problem has a solution provided that $\int_\G fdx=0$.
\end{Proposition}
\begin{Proof}
We shall proceed following the technique of \cite[Theorem 2.2]{nic} and of \cite[Theorem II.4.2]{ben}. To this end, it is expedient to introduce the following forms on $H^1(\G)$:
\begin{equation*}
a(u,v):= \sum_{j\in J}\int_{e_j}\left(\nu_j \pd_j u\pd_j v+ ub_j\pd_j v\right)dx,\qquad (u,v):= \sum_{j\in J}\int_{e_j}u_j v_jdx.
\end{equation*}
We observe that, for $s>0$ sufficiently large, the form
 \begin{equation}\label{a_s}
a_s(u,v):=a(u,v)+s(u,v)
\end{equation}
is coercive on $H^1(\G)$. Actually, the regularity of $b$ entails
\begin{eqnarray*}
a_s(u,u)&\geq& \sum_{j\in J}\int_{e_j}\left[\nu_j (\pd_j u)^2+s u_j^2-\|b\|_\infty |u_j| |\pd_j u|\right]dx\\
&\geq&  \sum_{j\in J}\int_{e_j}\left[\frac{\nu_j}2 (\pd_j u)^2+\left(s-\frac{\|b\|^2_\infty}{\nu_j} \right)u_j^2\right]dx.
\end{eqnarray*}
Fix $s>0$ such that $a_s$ is coercive on $H^1(\G)$. Invoking Lax-Milgram theorem, we obtain that for every $f\in L^2(\G)$, the problem
\begin{equation}\label{LMs}
a_s(u,v)=(f,v)\qquad \forall v\in H^1(\G)
\end{equation}
has exactly one solution $u=:G_s(f)$. In other words, there exists a map $G_s:L^2(\G)\to H^1(\G)$ such that $G_s(f)$ is the unique solution in $H^1(\G)$ of problem \eqref{LMs}. In fact we claim that
\[
G_s(f)\in H^2(\G)\qquad \forall f\in L^2(\G).
\]
Actually, let us observe that, on each edge $e_j$, the weak formulation of problem \eqref{fpCLA} is equivalent to the equality (in distributional sense)
\[
\nu_j \pd^2_jm=f+\pd_j(bm)
\]
where the right-hand side is in $L^2(e_j)$. Whence, $m\in H^2(e_j)$ for every $j\in J$ and our claim is completely proved.
Let us observe that the weak formulation is equivalent to
\[
(I-sG_s)(u)=G_s(f) \qquad \textrm{in }L^2(\G);
\]
indeed, the weak formulation can be written as: $a_s(u,v)=(f,v)+s(u,v)$;
hence $u=G_s(f)+s G_s(u)$.
We observe that, by the Rellich-Kondrachov theorem (see \cite{nic}), $G_s$ maps compactly $L^2(\G)$ into itself. By Fredholm alternative, the existence and the uniqueness of our problem are related to the properties of the operator $(I-sG_s^*)$ where $G_s^*$ is the adjoint operator of $G_s$.

Let us now calculate $G_s^*$. To this end, it is expedient to introduce the problem ($s>0$ is the same as before)
\begin{equation}\label{dual}
\left\{\begin{array}{ll}
sw-\nu\pd^2 w+b(x)\pd w=h &x\in \G\\
w_j(v_i)=w_k(v_i),\quad \sum_{j\in Inc_i} \nu_j\pd_j m(v_i)=0 & v_i\in\Gin, \, j,k\in Inc_i
\end{array}\right.
\end{equation}
whose weak formulation is
\[
\sum_{j\in J}\int_{e_j}\left(\nu_j \pd_j w\pd_j v+ b_j\pd_jw v+swv\right)dx=  \sum_{j\in J}\int_{e_j} h_j v dx\qquad \forall v\in H^1(\G).
\]
Arguing as before, we infer that there exists a compact map $\tilde G_s:L^2(\G)\to H^2(\G)$ such that $\tilde G_s(h)$ is the unique weak solution to problem~\eqref{dual}.

We claim that $G_s^*=\tilde G_s$. In order to prove this claim, it suffices to show that there holds $(G_s(f),h)=(f, \tilde G_s(h))$ for every $f,h\in L^2(\G)$. For $m:=G_s(f)$ and $z:=\tilde G_s(h)$, the regularizing effect of $G_s$ and $\tilde G_s$ ensures
\begin{eqnarray*}
(f, \tilde G_s(h))&=&\int_\G fz\, dx=\sum_j\int_{e_j}[s m-\nu_j \pd_j^2m-\pd(bm)]z\,dx\\
&=& \sum_j\int_{e_j}[smz +\nu_j \pd_jm\pd_jz +bm\pd_jz]\,dx=\int_\G mh\, dx=\int_\g G_s(f)h\,dx
\end{eqnarray*}
where the third equality is due to the Kirchhoff condition in~\eqref{fpCLA}, while the last two equalities are due respectively to the weak definition of~\eqref{dual} and to the definition of $m$. Hence, we accomplished the proof that $\tilde G_s$ is the dual of $G_s$.
Therefore, invoking Fredholm alternative, we infer that the dimension of the kernel of $(I-sG_s)$ is finite and it coincides with the one of the kernel of $(I-s\tilde G_s)$. in order to evaluate the latter, we observe that it is the space of the weak solution to~\eqref{dual} with $h=0$ and $s=0$.
For this problem, Theorem \cite[Theorem 2.1]{nic} ensures the maximum principle and in particular that every solution is constant. Therefore the kernel of $(I-sG_s)$ is one-dimensional and the first part of the statement is completely proved.\\
In conclusion, still Fredholm alternative guarantees
\[
Rg(I-sG_s)=Ker(I-sG_s^*)^\perp;
\]
in particular, problem~\eqref{fpCLA} has a solution provided that $f$ is orthogonal to a one-dimensional space. Using $\phi=1$ as test function  in~\eqref{weakCLA}, we obtain the desired compatibility condition for the solvability of the problem.
\end{Proof}
\begin{Lemma}\label{lemma:pdfp}
For every $\phi\in L^2(\G)$, the parabolic Cauchy problem
\begin{equation}\label{dualpara}
\left\{
\begin{array}{lll}
&\pd_t U- \nu \pd^2 U +b \pd U=0\qquad & \textrm{in }(0,+\infty)\times\G\\[4pt]
&\ds\sum_{j\in Inc_i} \nu_j\pd_j U(t,v_i)=0 & v_i\in\Gin,\, t\in(0,+\infty)\\
& U_j(t,v_i)=U_k(t,v_i) & v_i\in\Gin,\, j,k\in Inc_i,\, t\in(0,+\infty)\\
&U(0,x)=\psi(x) \qquad & \textrm{on }\G.
      \end{array}
      \right.
\end{equation}
 admits exactly one weak solution. Moreover, for $\phi\geq 0$, there holds:
\[
U(t,x)>0\qquad\forall (t,x)\in (0,+\infty)\times\G.
\]
\end{Lemma}
\begin{Proof}
The existence and uniqueness of the solution are established in \cite[Theorem 3.4]{nvb}.
For $\phi \geq 0$, the solution $U$ is strictly positive in $(0,+\infty)\times\G$ because the corresponding semigroup is holomorphic, positive and irreducible (see \cite{en} or \cite{ar}). We observe that the positivity of the semigroup is a straightforward consequence of the comparison principle. Moreover, the semigroup is holomorphic because the form~$a_s$ introduced in~\eqref{a_s} is coercive (as established in the proof of Proposition~\ref{ben4.2}). Finally, the irreducibility of the semigroup can be obtained following the same arguments of~\cite[Proposition 5.2]{kms}.
\end{Proof}
\begin{Remark}
Let us recall that estimates for the kernel function for problem \eqref{dualpara} can be obtained arguing as in \cite{r,c}.
\end{Remark}

\section*{Acknowledgment}
The authors are grateful to prof. Delio Mugnolo for several useful suggestions.
The second author has been partially supported by Gnampa-Indam and by Prat ``Traffic Flow on Networks: Analysis and Control'' by the University of Padova.



\begin{thebibliography}{99}

\bibitem{acc}
Y. Achdou, F. Camilli  and I. Capuzzo Dolcetta, Mean field games: numerical methods for the planning problem, SIAM J. Control Optim. 50 (2012), 77--109.

\bibitem{ar}
W. Arendt, \emph{Heat kernels}, Publication of the $9^{\textrm{th}}$ Internet Seminar 2005-06 of Ulm University, avaliable at www.uni-ulm.de/fileadmin/website\_uni\_ulm/mawi.inst.020/arendt/downloads/internetseminar.pdf

\bibitem{bf}
M. Bardi and E. Feleqi, Nonlinear elliptic systems and Mean Field Games, preprint avaliable at cvgmt.sns.it/paper/2655.

\bibitem{ben} A. Bensoussan, \emph{Perturbation methods in optimal control}, John Wiley and Sons, Chichester 1988.

\bibitem{bmp81}
L.Boccardo,   F. Murat and  J.P. Puel,  Existence de solutions faibles pour des \'equations elliptiques quasi-lin\'eaires \`a croissance quadratique. \emph{ Nonlinear partial differential equations and their applications.
 Coll\`ege de France Seminar, Vol. IV (Paris, 1981/1982)}, 19--73, Res. Notes in Math., 84, Pitman, Boston,  1983.

\bibitem{ccm}
F. Camilli, E. Carlini and C. Marchi, A model problem for Mean Field Games on networks,  Discrete Contin. Dyn. Syst. 35 (2015), no. 9, 4173--4192.

\bibitem{cms}
F. Camilli, C. Marchi and D. Schieborn, The vanishing viscosity limit for Hamilton-Jacobi equation on networks. J. Differential Equations 254 (2013), 4122--4143.

\bibitem{cgpt}
P. Cardaliaguet, J. Graber, A. Porretta and  D. Tonon, Second order mean field games with degenerate diffusion and local coupling, arXiv:1407.7024, 2014.

\bibitem{cgr}
P. Cardaliaguet and J. Graber, Mean field games systems of first order, Arxiv:1401.1789, 2014.

\bibitem{c}
C. Cattaneo, The spread of the potential on a weighted graph, Rend. Sem. Mat. Univ. Politec. Torino 57 (1999), no. 4, 221--229.

\bibitem{cg}
G.M. Coclite and M. Garavello, Vanishing viscosity for traffic on networks. SIAM J. Math. Anal. 42 (2010), no. 4, 1761--1783.

\bibitem{en}
K.-J. Engel and R. Nagel, \emph{A short course on operator semigroups}. Springer, New York, 2006.

\bibitem{fs}
M. Freidlin and  S. Sheu, Diffusion processes on graphs: stochastic differential equations, large deviation principle.  Probab. Theory Related Fields 116 (2000), no. 2, 181--220.

\bibitem{fw}
M. Freidlin and A. Wentzell, Diffusion processes on graphs and the averaging principle. Ann. Probab. 21 (1993), no. 4, 2215--2245.

\bibitem{gs}
D. Gomes and  J. Saude,  Mean field games - A brief survey.    Dyn. Games Appl. 4 (2014), no. 2, 110-154.

\bibitem{gms1}
D. Gomes, J. Mohr and R. Souza, Continuous time finite state mean field games. Appl. Math. Optim. 68 (2013), no. 1, 99--143.

\bibitem{gms2}
D. Gomes, J. Mohr and R. Souza, Discrete time, finite state space mean field games. J. Math. Pures Appl. (9) 93 (2010), no. 3, 308--328.

\bibitem{g.ppt}
O. Gu{\'e}ant, Existence and uniqueness result for mean field games with
congestion effect on graphs, Appl. Math. Optim. (to appear), http://arxiv.org/abs/1110.3442.

\bibitem{kms}
M. Kramar Fijavz, D. Mugnolo and  E.Sikolya,
 Variational and semigroup methods for waves and diffusion in networks. Appl. Math. Optim. 55 (2007), no. 2, 219--240.

\bibitem{hmc}
M. Huang, R.P. Malham\'e   and P.E. Caines,   Large population stochastic dynamic games: closed-loop McKean-Vlasov systems and the Nash certainty equivalence principle. Commun. Inf. Syst. 6 (2006), no. 3, 221--251.

\bibitem{ll}
J.-M. Lasry and P.-L. Lions,  Mean field games. Jpn. J. Math., 2 (2007), no. 1, 229--260.

\bibitem{lu}
 G.Lumer,  Espaces ramifi\'es  et diffusions sur les r\'eseaux topologiques.  C. R. Acad. Sci. Paris S\'er. A-B 291 (1980), no. 12, A627--A630.

\bibitem{mug}
D. Mugnolo, \emph{Semigroups methods for evolution equations on networks}, Understanding Complex Systems, Springer-Verlag,
Berlin, 2014.

\bibitem{nic}
S. Nicaise,  Elliptic operators on elementary ramified spaces. Integral Equations Operator Theory 11 (1988), no. 2, 230--257.

\bibitem{pb}
 Yu. V. Pokornyi and A. V. Borovskikh, Differential equations on networks (geometric graphs). J. Math. Sci. (N.Y.) 119 (2004), no. 6, 691--718.

\bibitem{p}
A. Porretta, On the Planning Problem for the Mean Field Games System, Dyn. Games Appl. (2014), no. 4, 231--256

\bibitem{r}
J.P. Roth, Le spectre du laplacien sur un graphe, in \emph{ Theorie du potentiel (Orsay, 1983)}, 521--539, Lecture Notes in Math., 1096, Springer, Berlin, 1984.

 \bibitem{vb}
J. von Below, Classical solvability of linear parabolic equations on networks. J. Differential Equations 72 (1988), no. 2, 316--337.

\bibitem{nvb}
J. von Below and S. Nicaise, Dynamical interface transition in ramified media with diffusion,  Comm. Partial Differential Equations 21 (1996), 255--279.


\end{thebibliography}
\end{document}